\documentclass[preprint,12pt]{elsarticle}
\usepackage{amssymb}
\usepackage{amsmath}
\usepackage{color}
\usepackage{amsthm}
\newtheorem{lm}{{\bf {Lemma}}}[section]
\newtheorem{thm}{{\bf {Theorem}}}[section]

\begin{document}
\begin{frontmatter}

\title{Ten limit cycles near a cubic homoclinic loop\\ with a nilpotent cusp}

\author[]{Yun Tian\corref{cor1}}
\ead{ytian22@shnu.edu.cn}

\author[]{Didi Ma}
%\ead{frank_gw@126.com}

\cortext[cor1]{Corresponding author. }
\address{Department of Mathematics, Shanghai Normal University, Shanghai, China}

\begin{abstract}
In this paper, we study the bifurcation of limit cycles near a homoclinic cuspidal loop
in a planar cubic near-Hamiltonian system by high-order Melnikov functions.
We present a method combining the algebraic structure of Abelian integrals and Picard-Fuchs equation
 for computing the corresponding asymptotic expansion of Melnikov functions near the cuspidal loop.
Using this system as an example, we show that planar cubic systems can have ten limit cycles
bifurcating near a cubic homoclinic loop.
\end{abstract}

\end{frontmatter}

\section{Introduction and main results }
One of restricted versions of the second part of Hilbert's 16th problem is to estimate the number of
limit cycles bifurcating near a homoclinic or heteroclinic loop in planar polynomial near-integrable systems of the form
\begin{equation}\label{a1}
\begin{split}
\dot x =&\, \mu^{-1}(x,y)H_y(x,y) + \varepsilon p(x,y,\varepsilon),\\
\dot y =&\, -\mu^{-1}(x,y)H_x(x,y) + \varepsilon q(x,y,\varepsilon),
\end{split}
\end{equation}
where $|\varepsilon|\ll1$, $H(x,y)$ is a smooth first integral of the unperturbed system \eqref{a1}$|_{\varepsilon=0}$
with the integrating factor $\mu(x,y)$, $p(x,y,\varepsilon)$ and $q(x,y,\varepsilon)$ are polynomials in $(x,y)$ and analytic in $\varepsilon$.
Suppose the integrable system \eqref{a1}$|_{\varepsilon=0}$ has a continuous family of periodic orbits given by
\[\Gamma_h\!:\, H(x,y)=h,\, h\in(\alpha, \beta),\]
with a homoclinic or heteroclinic loop $\Gamma_{\beta}\subset \{H(x,y)=\beta\}$.
 %as a part of its boundary.
Then the number of limit cycles produced near $\Gamma_{\beta}$ can be estimated
by the number of isolated zeros of the Melnikov function
\begin{equation*}
M(h) = \oint_{\Gamma_h} \mu(x,y)(q(x,y,0)\mathrm{d}x-p(x,y,0)\mathrm{d}y)%, \,\, h\in(\alpha, \beta),
\end{equation*}
for $0<\beta-h\ll1$.

%For limit cycle bifurcations near homoclinic or heteroclinic loops in near-integrable systems,
To determine the number of isolated zeros of $M(h)$ for $0<\beta-h\ll1$,
lots of results have been derived for the asymptotic expansion of $M(h)$
near $\Gamma_\beta$, for example see \cite{AZK2012,D1923,HYX2012,HZY2009,R1986,SH2019,XH2021,YYH2019,ZHX2008}.
Formulas were obtained for the first several coefficients of related asymptotic expansions of Melnikov functions
near different types of loops in \cite{HYT2008,HYX2012,HY1998,HZY2009}.
%Noticing that $H(3)\ge 13$, to find more limit cycles near the loops we need to get more coefficients.
Tian and Han \cite{TH2017} established an algorithm to compute more coefficients of the expansion
of Melnikov functions near homoclinic loops with studying bifurcations of small limit
cycles near an elementary center simultaneously.
This method was extended to some other cases
of heteroclinc loops \cite{GT2021}, cuspidal loops \cite{LH2020} and loops passing through
nilpotent saddles \cite{YYH2021,YYH2019}.

For the quadratic case of system \eqref{a1}, there are some known results obtained for limit cycle
bifurcations near a homoclinic or heteroclinic loop, for example see \cite{GI2015, H1997, HY2005,HYZ1999, HZ1999, HL1997,HI1994,I1996, XHX2021}
and references therein. It was proved that
the cyclicity of a homoclinic loop is two for quadratic Hamiltonian systems under quadratic perturbations
\cite{HYZ1999,HI1994,I1996}.
Gavrilov and Iliev \cite{GI2015} proved that there are at most three limit cycles bifurcating from
a heteroclinic loop with two saddles in quadratic near-Hamiltonian systems.
Xiong, Han and Xiao \cite{XHX2021} showed that three limit cycles can be produced near a triangle loop
in quadratic near-integrable systems.
For cubic systems, there are few results for finding lower bounds on the number of limit cycles near
homoclinic or heterclinic loops.
Examples were presented to show the existence of 5 limit cycles
near a homoclinic loop in \cite{TX2010}, and near a heteroclinic loop in \cite{GHTK2023}.

In this paper, we shall present an example of a planar cubic Hamiltonian system with a cuspidal homoclinic loop
from which ten limit cycles bifurcate under suitable cubic perturbations.
Consider the following cubic near-Hamiltonian system
\begin{equation}\label{sys}
\dot x = H_y(x,y)+\sum^{3}_{k=1}\varepsilon^k P_k(x,y),\quad
\dot y = - H_x(x,y)+\sum^{3}_{k=1}\varepsilon^k Q_k(x,y),
\end{equation}
where  $|\varepsilon|\ll 1$, $H(x,y)=\frac{1}{2}y^2 - \frac{1}{3}x^3 + \frac{1}{4}x^4$ and
\begin{equation}\label{b1}
P_k(x,y)=\sum_{i+j=1}^3 p_{ijk}x^iy^j,\quad
Q_k(x,y)=\sum_{i+j=1}^3 q_{ijk}x^iy^j
\end{equation}
with the coefficients $p_{ijk}$ and $q_{ijk}$ as free parameters.

The unperturbed system of \eqref{sys} has
an elementary center at the point $(1,0)$ and a nilpotent cusp at the origin.
There is a homoclinic loop $\Gamma_0$ defined by $H(x,y)=0$ passing
through the origin.  There exist two continuous families of periodic orbits
\begin{equation}\label{b1a}
\begin{split}
&\Gamma_h^+:~ H(x,y)=h,\quad h\in(0,+\infty), \\[0.5ex]
&\Gamma_h^-:~ H(x,y)=h,\quad h\in(-1/12,0), %0<\pm h\ll 1,
\end{split}
\end{equation}
which are separated by the cuspidal loop $\Gamma_0$.

Polynomial perturbations of the Hamiltonian system \eqref{sys}$|_{\varepsilon=0}$ were
investigated in \cite{LH2020,TX2010,WZ2020,ZZ1998}.
Zhao and Zhang \cite{ZZ1998} gave an upper bound $[\frac{7n}{2}]$ for the number of
zeros of the first order Melnikov functions for perturbations of degree $n$
on the two period annulus. %When $n=3$, the upper bound is $10$.
In \cite{TX2010} the existence of $5$ limit cycles near $\Gamma_0$
was proved for cubic perturbations by using the first order Melnikov
function. Here we shall study higher-order Melnikov functions
and find $10$ limit cycles near $\Gamma_0$ for suitable cubic perturbations.

Let $d^\pm(h,\varepsilon)$ be two displacement functions near $\Gamma_0$ for $\pm h>0$, respectively.
Then $d^\pm(h,\varepsilon)$ can be expanded as
\begin{equation}\label{b2}
d^\pm(h,\varepsilon)=\varepsilon M_1^\pm(h) + \varepsilon^2 M_2^\pm(h)
+\cdots+\varepsilon^k M_k^\pm(h)+\cdots,
\end{equation}
where
\begin{equation}\label{b2a}
M_1^\pm (h)=\oint_{\Gamma_h^\pm}Q_1(x,y)\mathrm{d}x-P_1(x,y)\mathrm{d}y.
\end{equation}
%The coefficients $M_k^{\pm}(h)$ of \eqref{b2} are called the $k$th-order Melnikov functions of system \eqref{sys}.
By studying  the first nonvanishing Melnikov functions in \eqref{b2}
we get the following two theorems.

\begin{thm}\label{thm1}
Let \eqref{b2} hold.  Then the following statements hold:

{\rm(I)} If $M_1^\pm(h)\not\equiv0$, system \eqref{sys}
can have 5 limit cycles bifurcating near the homoclinic loop $\Gamma_0$ with proper perturbations.
Furthermore, $M_1^\pm(h)\equiv0$ if and only if
\begin{equation}\label{b17}
p_{121}=-3q_{031},\quad q_{011}=-p_{101},\quad q_{111}=-2p_{201},\quad q_{211}=-3p_{301}.
\end{equation}

{\rm(II)} Assume $M_1^\pm(h)\equiv0$. If $M_2^\pm(h)\not\equiv0$, system \eqref{sys}
can have 8 limit cycles bifurcating near the homoclinic loop $\Gamma_0$ with proper perturbations.
Furthermore, $M_2^\pm(h)\equiv0$ if and only if one of the following conditions holds:
\begin{equation}\label{a4}
\begin{split}
&S_{1}:\left\{
\begin{aligned}
p_{111}=-2q_{021},\quad &p_{122}=-3q_{032},\quad q_{012}=-p_{102},\\
q_{112}=-2p_{202},\quad &q_{121}=-p_{211},\quad q_{212}=-3p_{302};
\end{aligned}
\right.\\[1ex]
&S_{2}:\left\{
\begin{aligned}
p_{122}=&\,\, p_{021}(p_{111}+2q_{021})-3q_{032},\quad q_{012} = -p_{102},\quad q_{031} = 0,\\
q_{112}=&\,\,p_{101}(p_{111}+2q_{021})-2p_{202},\quad q_{121} = -p_{211},\\
q_{212}=&\,\, (p_{201}+p_{301})(p_{111}+ 2q_{021})-3p_{302};
\end{aligned}
\right.\\[1ex]
&S_{3}:\left\{
\begin{aligned}
p_{021}=&\,0,\quad p_{122} = 2p_{301}(p_{211}+q_{121})-3q_{032},\quad q_{012} = -p_{102},\\
q_{031}=&\,0,\quad q_{112} = p_{101}(p_{111}+2q_{021})-2p_{202},\\
q_{212}=&\,2p_{101}(p_{211}+q_{121})+2(p_{201}+p_{301})(p_{211}+q_{021}+q_{121})\\
&+p_{111}(p_{201}+p_{301})-3p_{302}.
\end{aligned}
\right.
\end{split}
\end{equation}
\end{thm}

%In the next theorem, we show the existence of $10$ limit cycles near $\Gamma_0$ for system \eqref{sys}
%by studying $M_3^{\pm}(h)$.

\begin{thm}\label{thm3}
Let \eqref{b2} hold. If $M_1^\pm(h)=M_2^\pm(h)\equiv0,\,M^\pm_3(h)\not\equiv0$, system \eqref{sys}
can have 10 limit cycles produced near the homoclinic loop $\Gamma_0$ with proper perturbations.
\end{thm}

%For the bifurcation of limit cycles near $\Gamma_0$,
To prove Theorems \ref{thm1} and \ref{thm3},
we use the algebraic structure of Abelian integrals
\begin{equation*}
I_{i,j}^{\pm}(h)=\oint_{\Gamma_h^\pm}x^iy^j\mathrm{d}x
\end{equation*}
and Picard-Fuchs equation to
 study the asymptotic expansion of
 $M_k^{\pm}(h)$ near $h=0$, $k=1,2,3$.
%Formulas were only derived for the first five coefficients of the expansion
%of Melnikov functions near a homoclinic cuspidal loop.
%Using our method,
Recursive formulas are obtained for the coefficients of the expansion
of Melnikov functions near $\Gamma_0$
without investigating Hopf bifurcation near
the center simultaneously.

The paper is organized as follows: In Section \ref{sec2} we study the algebraic structure
of Abelian integrals $I_{i,j}^\pm(h)$ in Lemma \ref{lm1}
and obtain recursive formulas for the coefficients of the asymptotic expansion of $I_{i,1}^\pm(h)$, $i=0,1,2$, in Lemma \ref{lm2b}
by Picard-Fuchs equation.
Then we can get the asymptotic expansion of any linear combination of Abelian integrals $I_{i,j}^\pm(h)$ near $h=0$.
In Section \ref{sec3} we give the proofs of Theorems \ref{thm1} and \ref{thm3}, respectively.

\section{Preliminaries}\label{sec2}
We begin by studying Abelian integrals of system \eqref{sys}
\begin{equation}
I_{i,j}^{\pm}(h)=\oint_{\Gamma_h^\pm}x^iy^j\mathrm{d}x,
\nonumber
\end{equation}
where $\Gamma_h^\pm$ are given in \eqref{b1a}. %\!\subset\!\{\frac{1}{2}y^2 - \frac{1}{3}x^3 + \frac{1}{4}x^4=h\}$ for $0<\pm h\ll 1$.
%Note that $I^+_{i,j}(0)=I^-_{i,j}(0)$ for any integers $i\geq0$ and $j\geq-1$.
%Then for the simplicity of notation  we define functions $I_{i,j}(h)$ such that
%$I_{i,j}(h)=I^+_{i,j}(h)$ for $h\geq0$ and $I_{i,j}(h)=I^-_{i,j}(h)$ for $h\leq0$.
For the algebraic structure of Abelian integrals $I^\pm_{i,j}(h)$ we have the following lemma.

\begin{lm}\label{lm1}
For Abelian integrals $I^\pm_{i,j}(h)$ of system \eqref{sys}, the following identities hold:\vskip 1ex

{\rm(i)} $I^\pm_{i,j}(h)\equiv0$ for $j$ even;\vskip 1ex

{\rm(ii)} $(i+1)I^\pm_{i,j}(h)=j(I^\pm_{i+4,j-2}(h)-I^\pm_{i+3,j-2}(h))$ for $i\ge-1$ and $j\geq 1$;\vskip 1ex

{\rm(iii)} $I^\pm_{3,1}(h)=I^\pm_{2,1}(h)$ and
$I^\pm_{i,1}(h)=\frac{4i+6}{3i+9}I^\pm_{i-1,1}(h)+\frac{4i-12}{i+3}hI^\pm_{i-4,1}(h)\,\, \mbox{for}\,\, i\ge 4$.
\end{lm}

\begin{proof}
Because all the periodic orbits $\Gamma_h^{\pm}$ %:H= \frac{1}{2}y^2- \frac{1}{3}x^3+ \frac{1}{4}x^4=h$
are symmetric with respect to the $x$-axis, it is straightforward to get the statement $ (\rm{i})$.

Note that $ x^4/4=H(x,y)- y^2/2+ x^3/3$. For $i\geq-1$ and $ j\geq1$, we have
\begin{equation*}
\begin{split}
x^{i+4}y^{j-2}\mathrm{d}x&=x^{i+1}y^{j-2}\mathrm{d}(x^4/4)=x^{i+1}y^{j-2}\mathrm{d}\left(H- y^2/2+ x^3/3\right)\\
        &=x^{i+1}y^{j-2}\mathrm{d}H-x^{i+1}y^{j-1}\mathrm{d}y+x^{i+3}y^{j-2}\mathrm{d}x\\
        &=x^{i+1}y^{j-2}\mathrm{d}H-\mathrm{d}(x^{i+1}y^{j}/j)+\frac{i+1}{j}x^{i}y^{j}\mathrm{d}x+x^{i+3}y^{j-2}\mathrm{d}x.
\end{split}
\end{equation*}
Integrating its both sides along $\Gamma_h^\pm$ yields the statement $\rm(ii)$.

Taking $j=3$ and replacing $i$ by $i-4$ in $\rm(ii)$, for $i\ge 3$ we have
\begin{equation*}
\begin{split}
I^\pm_{i,1}(h)&=I^\pm_{i-1,1}(h)+\frac{i-3}{3}I^\pm_{i-4,3}(h)\\
        &=I^\pm_{i-1,1}(h)+\frac{i-3}{3}\oint_{\Gamma^\pm_h}\left(2h+\frac{2}{3}x^3-\frac{1}{2}x^4\right)x^{i-4}y\mathrm{d}x\\
       &=I^\pm_{i-1,1}(h)+\frac{2(i-3)}{3}hI^\pm_{i-4,1}(h)+\frac{2(i-3)}{9}I^\pm_{i-1,1}(h)-\frac{i-3}{6}I^\pm_{i,1}(h),
\end{split}
\end{equation*}
which yields the statement $\rm(iii)$. The proof is completed.
\end{proof}

Then by Lemma \ref{lm1} for any $i\ge0$ and $j\ge1$,
$I^\pm_{i,j}(h)$ can be expressed in the form
$$I^\pm_{i,j}(h)=p_1(h)I^\pm_{0,1}(h)+p_2(h)I^\pm_{1,1}(h)+p_3(h)I^\pm_{2,1}(h),$$
where $p_i(h)$, $i=1,2,3$, are polynomials in $h$.
Thus in order to expand $I^\pm_{i,j}(h)$ near $h=0$, we only need to study the asymptotic expansions of $I^\pm_{0,1}(h)$, $I^\pm_{1,1}(h)$
and $I^\pm_{2,1}(h)$.
To this end, we derive the Picard-Fuchs equation for  $I^\pm_{i,1}(h)$, $i=1,2,3$ in the next lemma.

\begin{lm}\label{lm2a}
Let $\mbox{\boldmath $X$}^\pm(h)=(I_{0,1}^\pm(h),\,I_{1,1}^\pm(h),\,I_{2,1}^\pm(h))^{T}$. Then
$\mbox{\boldmath $X$}^\pm(h)$ satisfies the Picard-Fuchs equation
\begin{equation}\label{b7}
P(h)\frac{\mathrm{d}}{\mathrm{d}h}\mbox{\boldmath $X$}^\pm(h)=(\mbox{\boldmath $A$}_1h+\mbox{\boldmath $A$}_0)\mbox{\boldmath $X$}^\pm(h),%\quad \pm h>0,
\end{equation}
where $P(h)=12h(12h+1)$, %{\boldmath $A$}$_0$ and {\boldmath $A$}$_1$ are two matrices given by
\begin{equation}\label{b8}
\mbox{\boldmath $A$}_0 =
\left( \begin{array}{ccc}
10 & 2 & -15 \\
0 & 14 & -15 \\
0 & 0 & 0
\end{array} \right), \quad
\mbox{\boldmath $A$}_1 =
\left( \begin{array}{ccc}
108 & 0 & 0 \\
-12 & 144 & 0 \\
-12 & -24 & 180
\end{array} \right).
\end{equation}
\end{lm}
%\begin{equation}\label{b3}
%12h(12h+1)\left(\!\begin{array}{c}
%I'_{0,1}\\[1ex] I'_{1,1}\\[1ex] I'_{2,1}
%\end{array}\!\right)=
%\left(\!\begin{array}{ccc}
%108h\!+\!10 & 2 & -15\\[1ex]
%-12h & 144h\!+\!14 & -15\\[1ex]
%-12h & -24h & 180h
%\end{array}\!\right)
%\left(\!\begin{array}{c}
%I_{0,1}\\[1ex] I_{1,1}\\[1ex] I_{2,1}
%\end{array}\!\right).
%\end{equation}
%\end{lm}

\begin{proof}
 Noting $y^2=2h+\frac{2}{3}x^3-\frac{1}{2}x^4 $ along $\Gamma_h^\pm$, we have
\begin{equation}\label{b4}
I^\pm_{i,1}(h)=\oint_{\Gamma_h^\pm}\frac{x^iy^2}{y}\,\mathrm{d}x=2hI^\pm_{i,-1}(h)+\frac{2}{3}I^\pm_{i+3,-1}(h)-\frac{1}{2}I^\pm_{i+4,-1}(h).
\end{equation}
By (ii)$|_{j=1}$ of Lemma \ref{lm1}, for $i\ge -1$ we get
\begin{equation}\label{b5}
I^\pm_{i+4,-1}(h)=(i+1)I^\pm_{i,1}(h)+I^\pm_{i+3,-1}(h).
\end{equation}
Then substituting \eqref{b5} into \eqref{b4} yields
\begin{equation}\label{b5a}
(3i+9)I^\pm_{i,1}(h)=12hI^\pm_{i,-1}(h)+I^\pm_{i+3,-1}(h) .
\end{equation}
When $i=-1,0,1$, by \eqref{b5} we can further get
\begin{equation}\label{b6}
\begin{split}
I^\pm_{3,-1}(h)=&\,I^\pm_{2,-1}(h),\\[1ex]
I^\pm_{4,-1}(h)=&\,I^\pm_{0,1}(h)+I^\pm_{3,-1}(h)=I^\pm_{0,1}(h)+I^\pm_{2,-1}(h),\\[1ex]
I^\pm_{5,-1}(h)=&\,2I^\pm_{1,1}(h)+I^\pm_{4,-1}(h)=2I^\pm_{1,1}(h)+I^\pm_{0,1}(h)+I^\pm_{2,-1}(h).
\end{split}
\end{equation}
Taking $i=0,1,2$ for \eqref{b5a} and then removing $I^\pm_{3,-1}(h)$, $I^\pm_{4,-1}(h)$ and $I^\pm_{5,-1}(h)$ by \eqref{b6},
 we obtain
\begin{equation*}
\begin{split}
9I^\pm_{0,1}(h)=&\,12hI^\pm_{0,-1}(h)+I^\pm_{2,-1}(h),\\
-I^\pm_{0,1}(h)+12I^\pm_{1,1}(h)=&\,12hI^\pm_{1,-1}(h)+I^\pm_{2,-1}(h),\\
-I^\pm_{0,1}(h)-2I^\pm_{1,1}(h)+15I_{2,1}(h)=&\,(12h+1)I^\pm_{2,-1}(h),\\
\end{split}
\end{equation*}
which yields \eqref{b7} because $\frac{\mathrm{d}}{\mathrm{d}h}I^\pm_{i,1}(h)=I^\pm_{i,-1}(h)$.
The proof is completed.
\end{proof}

Next, we shall use the Picard-Fuchs equation \eqref{b7} to compute the asymptotic expansions of $I_{0,1}^\pm(h),\,I_{1,1}^\pm(h)$ and $I_{2,1}^\pm(h)$
near $h=0$. % respectively.
By \cite{HZY2009}, $\mbox{\boldmath $X$}^\pm(h)$ can be expanded in the following form
\begin{equation}\label{b9}
\mbox{\boldmath $X$}^\pm(h)= \sum_{j=0}^{\infty} \mbox{\boldmath $a$}_j^\pm h^j+ \sum_{j=0}^{\infty}(\mbox{\boldmath $b$}_j^\pm {\lvert h \rvert}^\frac{5}{6}+\mbox{\boldmath $c$}_j^\pm {\lvert h \rvert}^\frac{7}{6}){\lvert h \rvert}^j,\quad 0\leq \pm h \ll1,
\end{equation}
where $\mbox{\boldmath $a$}_j^\pm$, $\mbox{\boldmath $b$}_j^\pm$ and $\mbox{\boldmath $c$}_j^\pm$ are vector coefficients. We have the next lemma.

\begin{lm}\label{lm2b}
Let $\mbox{\boldmath $X$}^\pm(h)=(I_{0,1}^\pm(h), I_{1,1}^\pm(h), I_{2,1}^\pm(h))^{T}$.
Then for the  vector coefficients in \eqref{b9} we have
\begin{equation}\label{b9a}
\begin{split}
&\mbox{\boldmath $A$}_0\,\mbox{\boldmath $a$}^\pm_0=0,\,\qquad
(\mbox{\boldmath $A$}_0-10\,\mbox{\boldmath $E$}_3)\mbox{\boldmath $b$}^\pm_0=0,\,\qquad
(\mbox{\boldmath $A$}_0-14\,\mbox{\boldmath $E$}_3)\mbox{\boldmath $c$}^\pm_0=0, \\[1ex]
&(\mbox{\boldmath $A$}_0-12j\,\mbox{\boldmath $E$}_3)\,\mbox{\boldmath $a$}^\pm_j=-(\mbox{\boldmath $A$}_1-144(j-1)\,\mbox{\boldmath $E$}_3)\,\mbox{\boldmath $a$}^\pm_{j-1},\ j\geq1;\\[1ex]
&(\mbox{\boldmath $A$}_0-(12j+10)\,\mbox{\boldmath $E$}_3)\,\mbox{\boldmath $b$}^\pm_j=\mp(\mbox{\boldmath $A$}_1-(144j-24)\,\mbox{\boldmath $E$}_3)\,\mbox{\boldmath $b$}^\pm_{j-1},\ j\geq1;\\[1ex]
&(\mbox{\boldmath $A$}_0-(12j+14)\,\mbox{\boldmath $E$}_3)\,\mbox{\boldmath $c$}^\pm_j=\mp(\mbox{\boldmath $A$}_1-(144j+24)\,\mbox{\boldmath $E$}_3)\,\mbox{\boldmath $c$}^\pm_{j-1},\ j\geq1.
\end{split}
\end{equation}
where $\mbox{\boldmath $A$}_0$ and $\mbox{\boldmath $A$}_1$ are matrices given in \eqref{b8}, and $\mbox{\boldmath $E$}_3$ is the $3\times3$ identity matrix.
\end{lm}

We shall use the Picard-Fuchs equation \eqref{b7} to prove Lemma \ref{lm2b}.
By \eqref{b9} we notice that the derivative $\frac{\mathrm{d}}{\mathrm{d}h}\mbox{\boldmath $X$}^\pm(h)$ is not defined at $h=0$
when $\mbox{\boldmath $b$}^\pm_0\neq0$. However,
the left-hand side of \eqref{b7} can be rewritten as
\begin{equation}\label{b11}
P(h)\frac{\mathrm{d}}{\mathrm{d}h}\mbox{\boldmath $X$}^\pm(h)
=\frac{\mathrm{d}}{\mathrm{d}h}(P(h)\mbox{\boldmath $X$}^\pm(h))-P'(h)\mbox{\boldmath $X$}^\pm(h),\,
\end{equation}
which can be expanded as a convergent series near $h=0$ by \eqref{b9}.
Then $P(h)\frac{\mathrm{d}}{\mathrm{d}h}\mbox{\boldmath $X$}^\pm(h)$ can be well-defined at $h=0$,
which implies that we can use \eqref{b7} and \eqref{b11}
to study the asymptotic expansion \eqref{b9} for $\mbox{\boldmath $X$}^\pm(h)$.

\begin{proof}
Here we only present the proof for $\mbox{\boldmath $X$}^+(h)$. The case for $\mbox{\boldmath $X$}^-(h)$ can be similarly proved.
By \eqref{b9} we have
\begin{equation*}
\begin{split}
P(h)\mbox{\boldmath $X$}^{+}(h)=&\ 12\mbox{\boldmath $a$}^{+}_0h+\sum_{j=2}^{\infty}(12\mbox{\boldmath $a$}^+_{j-1}+144\mbox{\boldmath $a$}^+_{j-2})h^j+12\mbox{\boldmath $b$}^+_0h^{\frac{11}{6}}+12\mbox{\boldmath $c$}^+_0h^{\frac{13}{6}}\\
&+\sum_{j=2}^{\infty}\left[(12\mbox{\boldmath $b$}^+_{j-1}+144\mbox{\boldmath $b$}^+_{j-2})h^{\frac{5}{6}}+(12\mbox{\boldmath $c$}^+_{j-1}+144\mbox{\boldmath $c$}^+_{j-2})h^{\frac{7}{6}}\right]h^j, \\
P'(h)\mbox{\boldmath $X$}^{+}(h)=&\ 12\mbox{\boldmath $a$}^{+}_0+\sum_{j=1}^{\infty}(12\mbox{\boldmath $a$}^+_{j}+288\mbox{\boldmath $a$}^+_{j-1})h^j+12\mbox{\boldmath $b$}^+_0h^{\frac{5}{6}}+12\mbox{\boldmath $c$}^+_0h^{\frac{7}{6}}\\
&+\sum_{j=1}^{\infty}\left[(12\mbox{\boldmath $b$}^+_{j}+288\mbox{\boldmath $b$}^+_{j-1})h^{\frac{5}{6}}+(12\mbox{\boldmath $c$}^+_{j}+288\mbox{\boldmath $c$}^+_{j-1})h^{\frac{7}{6}}\right]h^j.
\end{split}
\end{equation*}
Then by \eqref{b11}, the left-hand side of \eqref{b7} can be expanded as %in the following form
\begin{equation}\label{b12}
\begin{split}
P(h)\frac{\mathrm{d}}{\mathrm{d}h}\mbox{\boldmath $X$}^+(h)=&\,\sum_{j=1}^{\infty}[12j\mbox{\boldmath $a$}^+_{j}+144(j-1)\mbox{\boldmath $a$}^+_{j-1}]h^j+10\mbox{\boldmath $b$}^+_0h^{\frac{5}{6}}+14\mbox{\boldmath $c$}^+_{0}h^{\frac{7}{6}}\\
&+\,\sum_{j=1}^{\infty}[(12j+10)\mbox{\boldmath $b$}^+_{j}+(144j-24)\mbox{\boldmath $b$}^+_{j-1}]h^{j+\frac{5}{6}}\\
&+\,\sum_{j=1}^{\infty}[(12j+14)\mbox{\boldmath $c$}^+_{j}+(144j+24)\mbox{\boldmath $c$}^+_{j-1}]h^{j+\frac{7}{6}}.
\end{split}
\end{equation}
On the other hand, substituting \eqref{b9} into the right-hand side of \eqref{b7} yields
\begin{equation}\label{b13}
\begin{split}
&(\mbox{\boldmath $A$}_1h+\mbox{\boldmath $A$}_0)\mbox{\boldmath $X$}^+(h)
=\mbox{\boldmath $A$}_0\mbox{\boldmath $a$}^+_0+\sum_{j=1}^{\infty}(\mbox{\boldmath $A$}_0\mbox{\boldmath $a$}^+_j+\mbox{\boldmath $A$}_1\mbox{\boldmath $a$}^+_{j-1})h^j+\mbox{\boldmath $A$}_0\mbox{\boldmath $b$}^+_0h^{\frac{5}{6}}
\\ &\quad\,\,\,\,
+\!\mbox{\boldmath $A$}_0\mbox{\boldmath $c$}^+_0h^{\frac{7}{6}}
\!+\!\sum_{j=1}^{\infty}\left[(\mbox{\boldmath $A$}_0\mbox{\boldmath $b$}^+_{j}\!+\!\mbox{\boldmath $A$}_1\mbox{\boldmath $b$}^+_{j-1})h^{j+\frac{5}{6}}
\!+\!(\mbox{\boldmath $A$}_0\mbox{\boldmath $c$}^+_{j}\!+\!\mbox{\boldmath $A$}_1\mbox{\boldmath $c$}^+_{j-1})h^{j+\frac{7}{6}}\right].
\end{split}
\end{equation}
Comparing the coefficients in \eqref{b12} and \eqref{b13}, we can get \eqref{b9a} for $\mbox{\boldmath $a$}^+_j,\,\mbox{\boldmath $b$}^+_j$ and $\mbox{\boldmath $c$}^+_j,\,j\geq0$. The proof is completed.
\end{proof}

Note that by \eqref{b8} for any $j\geq1$ we have
\begin{equation*}
\mbox{det}(\mbox{\boldmath $A$}_0-(12j+\lambda)\mbox{\boldmath $E$}_3)\neq0,\quad\lambda\in\left\{0,10,14\right\}.
\end{equation*}
Then by the last three equations of \eqref{b9a}, for any $j\geq1$ the coefficients $\mbox{\boldmath $a$}^\pm_j,\,\mbox{\boldmath $b$}^\pm_j$ and $\mbox{\boldmath $c$}^\pm_j$
are uniquely determined by $\mbox{\boldmath $a$}^\pm_{j-1},\,\mbox{\boldmath $b$}^\pm_{j-1}$ and $\mbox{\boldmath $c$}^\pm_{j-1}$, respectively.
To be precise, for $j\ge1$ we have
\begin{equation}\label{b13b}
\begin{split}
\mbox{\boldmath $a$}^\pm_j=&-(\mbox{\boldmath $A$}_0-12j\,\mbox{\boldmath $E$}_3)^{-1}(\mbox{\boldmath $A$}_1-144(j-1)\,\mbox{\boldmath $E$}_3)\,\mbox{\boldmath $a$}^\pm_{j-1}, \\[1ex]
\mbox{\boldmath $b$}^\pm_j=&\mp(\mbox{\boldmath $A$}_0-(12j+10)\,\mbox{\boldmath $E$}_3)^{-1}(\mbox{\boldmath $A$}_1-(144j-24)\,\mbox{\boldmath $E$}_3)\,\mbox{\boldmath $b$}^\pm_{j-1},\\[1ex]
\mbox{\boldmath $c$}^\pm_j=&\mp(\mbox{\boldmath $A$}_0-(12j+14)\,\mbox{\boldmath $E$}_3)^{-1}(\mbox{\boldmath $A$}_1-(144j+24)\,\mbox{\boldmath $E$}_3)\,\mbox{\boldmath $c$}^\pm_{j-1}.
\end{split}
\end{equation}
Then when we have the values of $\mbox{\boldmath $a$}^\pm_0$, $\mbox{\boldmath $b$}^\pm_0$ and $\mbox{\boldmath $c$}^\pm_0$,
we can get all the remaining coefficients in \eqref{b9} by using \eqref{b13b}.

Suppose that $\mbox{\boldmath $a$}^\pm_0$, $\mbox{\boldmath $b$}^\pm_0$ and $\mbox{\boldmath $c$}^\pm_0$ are given by
\[\mbox{\boldmath $a$}^\pm_0=(a^\pm_{10},\,a^\pm_{20},\,a^\pm_{30})^T,\,\,\,
\mbox{\boldmath $b$}^\pm_0=(b^\pm_{10},\,b^\pm_{20},\, b^\pm_{30})^T,\,\,\,
\mbox{\boldmath $c$}^\pm_0=(c^\pm_{10},\,c^\pm_{20},\,c^\pm_{30})^T.\]
By the first three equations of \eqref{b9a}, we can get
\begin{equation}\label{b13c}
a^\pm_{20}=\frac{5}{6}a^\pm_{10},\quad
a^\pm_{30}=\frac{7}{9}a^\pm_{10},\quad
b^\pm_{20}=b^\pm_{30}=0,\quad
c^\pm_{20}=2c^\pm_{10},\quad c^\pm_{30}=0.
\end{equation}
Using the method presented in \cite{HZY2009}, we obtain
\begin{equation}\label{b13d}
\begin{split}
a^\pm_{10}=&\oint_{\Gamma_0}y\,\mathrm{d}x=2\int^{\frac{4}{3}}_0x^{\frac{3}{2}}\sqrt{\frac{2}{3}
-\frac{1}{2}x}\,\mathrm{d}x=\, \frac{4\sqrt{2}}{27}\pi,\\
b^\pm_{10}=&\,\widetilde{r}_{00}B^\pm_{00}=2\sqrt{2}\times(-\frac 1 3)^{-\frac{1}{3}}B^\pm_{00}=-2\,\sqrt[6]{72}B^\pm_{00},\\
%b^-_{10}=&\,\widetilde{r}_{00}{\color{red}B^-_{00}}=-2\,\sqrt[6]{72}{\color{red}B^-_{00}},\\
c^\pm_{10}=&\,\widetilde{r}_{10}B^\pm_{10}=-\frac{4}{3}\sqrt{2}\times(-\frac{1}{3})^{-\frac{5}{3}}\times\frac{1}{4}B^\pm_{10}
=\sqrt[6]{648}B^\pm_{10},
%,\\ c^-_{10}=&\,\widetilde{r}_{10}{\color{red}B^-_{10}}=\sqrt[6]{648}{\color{red}B^-_{10}}.
\end{split}
\end{equation}
where
\begin{equation}\label{b13e}
\begin{split}
B^+_{00}=&-\frac{3}{5}\int_{-\infty}^{1}\frac{\mathrm{d}x}{\sqrt{1-x^3}}<0,\qquad B^-_{00}=\frac{3}{5}\int_{0}^{1}\frac{\mathrm{d}x}{\sqrt{x(1-x^3)}}>0,\\ B^+_{10}=&\frac{3}{7}\left(\int^{-1}_1\frac{x\mathrm{d}x}{\sqrt{1-x^3}}-\int^1_0\frac{x^{\frac{3}{2}}\mathrm{d}x}{\sqrt{1+x^3}+1+x^3}-2\right)<0,\\
B^-_{10}=&-\frac{3}{7}\left(\int_{0}^{1}\frac{x^{\frac{3}{2}}\mathrm{d}x}{\sqrt{1-x^3}+1-x^3}-2\right)>0.
\end{split}
\end{equation}
Then by \eqref{b13b}, \eqref{b13c}
and \eqref{b13d} we can derive the following lemma.

\begin{lm}\label{lm5}
Abelian integrals $I^\pm_{0,1}(h)$, $I^\pm_{1,1}(h)$ and
$I^\pm_{2,1}(h)$ have the following asymptotic expansions near $h=0$:
\begin{equation}
\begin{split}\label{b13a}
I^\pm_{0,1}(h)=&\,\frac{4\sqrt{2}}{27}\pi-2 b^\pm_{0}|h|^{\frac{5}{6}} + b^\pm_{1}|h|^{\frac{7}{6}}
\,\pm\, \frac{35}{44}b^\pm_{0}|h|^{\frac{11}{6}} \,\mp\, \frac{385}{208}b^\pm_{1}|h|^{\frac{13}{6}}+\cdots,\\
I^\pm_{1,1}(h)=&\,\frac{10\sqrt{2}}{81}\pi + \sqrt{8}\pi h + 2b^\pm_{1}|h|^{\frac{7}{6}}
\,\pm\, \frac{21}{22}b^\pm_{0}|h|^{\frac{11}{6}} \,\mp\, \frac{55}{26}b^\pm_{1}|h|^{\frac{13}{6}}+\cdots,\\
I^\pm_{2,1}(h)=&\,\frac{28\sqrt{2}}{243}\pi+\frac{4\sqrt{2}}{3}\pi h \,\pm\, \frac{12}{11}b^\pm_{0}|h|^{\frac{11}{6}}
\,\mp\, \frac{30}{13}b^\pm_{1}|h|^{\frac{13}{6}}+\cdots,
\end{split}
\end{equation}
where $b^\pm_0=\sqrt[6]{72}\,B^\pm_{00}$ and $b^\pm_1=\sqrt[6]{648}\,B^\pm_{10}$
with $B^\pm_{00}$ and $B^\pm_{10}$ given in \eqref{b13e}.
\end{lm}

Then for any Melnikov functions $M^\pm(h,\mbox{\boldmath $\delta$})$ of system \eqref{sys} given by
% can be expressed as %linear combinations of Abelian integrals $I_{i,j}^\pm(h)$, $i\le i+j \le n$.
\begin{equation*}\label{c14}
M^\pm(h,\mbox{\boldmath $\delta$})=\sum_{i=0}^n \sum_{j=0}^{n-i} a_{ij}(\mbox{\boldmath $\delta$})I^\pm_{i,j}(h),\quad n\in \mathbb{N}^+,
\end{equation*}
we can use Lemmas \ref{lm1} and \ref{lm5} to compute the coefficients of the following asymptotic expansions
%then $M^\pm(h,\mbox{\boldmath $\delta$})$ can be expanded in the form
\begin{equation}\label{c15}
\begin{split}
M^+(h,\mbox{\boldmath $\delta$})
&=c_0(\mbox{\boldmath $\delta$})\! +\! \sum_{j\ge0}
\left(c_{3j+1}(\mbox{\boldmath $\delta$})h^{j+\frac{5}{6}}\! +\! c_{3j+2}(\mbox{\boldmath $\delta$})h^{j+1}
\! +\! c_{3j+3}(\mbox{\boldmath $\delta$})h^{j+\frac{7}{6}}\right),\,\\
M^-(h,\mbox{\boldmath $\delta$})
&=\tilde c_0(\mbox{\boldmath $\delta$})\! +\! \sum_{j\ge0}
\left(\tilde c_{3j+1}(\mbox{\boldmath $\delta$})|h|^{j+\frac{5}{6}}\! +\! \tilde c_{3j+2}(\mbox{\boldmath $\delta$})h^{j+1}
 \!+\!\tilde  c_{3j+3}(\mbox{\boldmath $\delta$})|h|^{j+\frac{7}{6}}\right),\\
\end{split}
\end{equation}
where $\mbox{\boldmath $\delta$}\in \mathbb{R}^m$ is a vector consisting of all the free parameters in system \eqref{sys}.
Note that by \eqref{b13b}, \eqref{b13c} and \eqref{b13d} we can get
\begin{equation*}
\mbox{\boldmath $a$}^+_j=\mbox{\boldmath $a$}_j^-,\quad
\mbox{\boldmath $b$}^+_j=(-1)^{j+1}\rho_1\, \mbox{\boldmath $b$}_j^-,\quad
\mbox{\boldmath $c$}^+_j=(-1)^{j+1} \rho_3\, \mbox{\boldmath $c$}_j^-,\quad j\ge 0,
\end{equation*}
where $\rho_1=-\frac{B_{00}^+}{B_{00}^-}>0$ and $\rho_3=-\frac{B_{10}^+}{B_{10}^-}>0$.
Then for \eqref{c15} using Lemma \ref{lm1}  we can further obtain
\begin{equation}\label{c16}
\begin{split}
&c_0=\tilde c_0, \quad c_{3j+2}=\tilde c_{3j+2}, \quad
 c_{3j+k}=(-1)^{j+1}\rho_k\,\tilde c_{3j+k},\quad
% c_{3j+3}=(-1)^{j+1}\tilde \rho \,\tilde c_{3j+3},
 k=1,3,\quad j\ge 0.
 \end{split}
\end{equation}

The identities in \eqref{c16} show the relation between the coefficients $c_k$ and $\tilde c_k$, $k\ge0$,
which can be also obtained by \cite{YH2023}.
In \cite{YH2023} Yang and Han presented the relation between the coefficients of the asymptotic expansions
of $M^\pm(h)$ near a cuspidal loop defined by an analytic Hamiltonian.
Based on this result, bifurcations of limit cycles near a cuspidal loop can be investigated
by only using the coefficients of the asymptotic expansion of $M^+(h)$.
The next lemma is derived from \cite{YH2023} for system \eqref{sys}.

\begin{lm}\label{thm4}
Assume that $M^\pm(h,\mbox{\boldmath $\delta$})$ are the first nonvanishing Melnikov functions in \eqref{b2}
with \eqref{c15} holding. If there exist $\mbox{\boldmath $\delta$}_0\in \mathbb{R}^m$ and
$3n+1\le k \le 3(n+1)$ such that
\begin{equation*}
\begin{split}
&c_0(\mbox{\boldmath $\delta$}_0)=c_1(\mbox{\boldmath $\delta$}_0)=\cdots=c_{k-1}(\mbox{\boldmath $\delta$}_0)=0,\quad
c_k(\mbox{\boldmath $\delta$}_0)\neq0,\\
&{\mathrm {rank}} \frac{\partial(c_0,c_1,\cdots,c_{k-1})}{\partial \mbox{\boldmath $\delta$}}\bigg|_{\mbox{\boldmath $\delta$}=\mbox{\boldmath $\delta$}_0}
=k,
\end{split}
\end{equation*}
then system \eqref{sys} can have $2k-(n+1)$ limit cycles near $\Gamma_0$ with $k$ (or $k-(n+1)$) limit cycles inside $\Gamma_0$
and $k-(n+1)$ (or $k$) limit cycles outside $\Gamma_0$ for some $(\varepsilon,\mbox{\boldmath $\delta$})$ near $(0,\mbox{\boldmath $\delta$}_0)$.
\end{lm}

\section{Proof of Theorems \ref{thm1} and \ref{thm3}}\label{sec3}

In this section, we shall prove Theorems \ref{thm1} and \ref{thm3}
by using Lemmas \ref{lm1}, \ref{lm5} and \ref{thm4}.

%\begin{proof}
\noindent{\bf Proof of Theorem \ref{thm1}.}
Note that Melnikov functions $M_1^\pm(h)$ in \eqref{b2a} can be written in the form
\begin{equation}\label{b2b}
M^\pm_1(h)= \oint_{\Gamma^\pm_h} \left(Q_1(x,y) + \int \frac{\partial}{\partial x} P_1(x,y)\mathrm{d}y\right)\mathrm{d}x.
\end{equation}
Then by \eqref{b1} and \eqref{b2b} we have
\begin{equation}\label{b14a}
M^\pm_1(h)=A_1I^\pm_{0,1}(h)+A_2I^\pm_{1,1}(h)+A_3I^\pm_{2,1}(h)+A_4I^\pm_{0,3}(h),
\end{equation}
where
\begin{equation}\label{b14b}
\begin{split}
&A_1=p_{101}+q_{011},\ A_2=2p_{201}+q_{111},\\
& A_3=3p_{301}+q_{211},\ A_4=\frac{1}{3}p_{121}+q_{031}.
\end{split}
\end{equation}

By Lemma \ref{thm4} we only need to study the coefficients of the asymptotic expansion
of $M^+_1(h)$ for $0\le h\ll 1$.
By the statements (ii) and (iii) of Lemma \ref{lm1}, we get
\begin{equation}\label{b14}
\begin{split}
I^+_{4,1}(h)=\,&\frac{4}{7}hI^+_{0,1}(h)+\frac{22}{21}I^+_{3,1}(h)=\frac{2}{21}(6hI^+_{0,1}(h)+11I^+_{2,1}(h)),\\
I^+_{0,3}(h)=\,&3I^+_{4,1}(h)-3I^+_{3,1}(h)=\frac{1}{7}(12hI^+_{0,1}(h)+I^+_{2,1}(h)).
\end{split}
\end{equation}
Substituting the second identity of \eqref{b14} into \eqref{b14a} yields
\begin{equation*}\label{b15}
M^+_1(h)=(A_1+\frac{12}{7}A_4h)\,I^+_{0,1}(h)+A_2\,I^+_{1,1}(h)+(A_3+\frac{1}{7}A_4)\,I^+_{2,1}(h).
\end{equation*}
Then using \eqref{b13a} we further obtain
\begin{equation}\label{b15a}
M^+_1(h)=c_0+b^+_{0}c_1h^{\frac{5}{6}}+c_2h+b^+_{1}c_3h^{\frac{7}{6}}+\cdots,\quad 0<h\ll 1,
\end{equation}
where
\begin{equation}\label{b16}
\begin{split}
c_0=&\,\frac{2\sqrt{2}\pi}{243}(18A_1+15A_2+14A_3+2A_4),
\quad c_1=-2A_1,\\[1ex]
c_2=&\,\frac{2\sqrt{2}\pi}{9}(9A_2+6A_3+2A_4),
\quad\quad c_3=A_1+2A_2.
\end{split}
\end{equation}
Then by \eqref{b14b} and \eqref{b16} solving $c_0=c_1=c_2=0$ in $q_{011},\,q_{111}$ and $q_{211}$ yields
\begin{equation*}
\begin{split}
& q_{011}=-p_{101},\quad q_{111}=-2p_{201}-\frac{4}{27}(p_{121}+3q_{031}),\\
& q_{211}=-3p_{301}+\frac{1}{9}(p_{121}+3q_{031}),
\end{split}
\end{equation*}
and $c_3=-\frac{8}{27}(p_{121}+3q_{031})$.
Because $c_0,\;c_1$ and $c_2$ are linear in $q_{011},\;q_{111}$ and $q_{211}$, by Lemma \ref{thm4}
system \eqref{sys} can have $5$ limit cycles near $\Gamma_0$ with proper perturbations with $p_{121}+3q_{031}\neq0$.

By \eqref{b16} it is easy to get
\begin{equation*}
c_j=0,\;j=0,1,2,3\ \Longleftrightarrow\ A_j=0,\;j=1,2,3,4.
\end{equation*}
Then by \eqref{b14a} and \eqref{b15a} we obtain $M^\pm_1(h)\equiv0$ if and only if $A_1=A_2=A_3=A_4=0$, which yields \eqref{b17}.

Next, we assume $M_1^\pm(h)\equiv0$, and then study the asymptotic expansion of $M_2^+(h)$ near $h=0$.
With \eqref{b17} holding, we get
\begin{equation*}
Q_1\mathrm{d}x-P_1\mathrm{d}y=r_1\mathrm{d}H+\mathrm{d}{R_1},
\end{equation*}
where
\begin{equation}\label{b17a}
r_1(x,y)=-(p_{111}+2q_{021})x-(p_{211}+q_{121})x^2.
\end{equation}
Then using Fran\c coise's algorithm \cite{JP1996}, we have
\begin{equation*}
M^\pm_2(h)=\oint_{\Gamma_h^\pm}Q\mathrm{d}x-P\mathrm{d}y,
\end{equation*}
where
\begin{equation*}
\begin{split}
Q(x,y)=Q_2(x,y)+r_1(x,y)Q_1(x,y),\\
P(x,y)=P_2(x,y)+r_1(x,y)P_1(x,y).
\end{split}
\end{equation*}

Then by the statement $\rm(i)$ of Lemma \ref{lm1}, $M^\pm_{2}(h)$ can be simplified into the form
\begin{equation}\label{b18}
\begin{split}
M_2^\pm(h)=& \oint_{\Gamma_h^\pm}\left(Q+\int \frac{\partial}{\partial x}P\mathrm{d}y \right)\mathrm{d}x\\
=&\sum_{j=1}^{5}B_jI^\pm_{j-1,1}(h)+\sum_{j=6}^{8}B_jI^\pm_{j-6,3}(h),
\end{split}
\end{equation}
where
\begin{equation}\label{b18a}
\begin{split}
B_1 =&\   p_{102}+q_{012},\qquad B_2 = 2p_{202}+q_{112}-(p_{111}+2q_{021})p_{101},\\
B_3 =&\   3p_{302}+q_{212}-(p_{111}+2q_{021})p_{201}-(2p_{211}+2q_{121})p_{101},\\
B_4 =& -(p_{111}+2q_{021})p_{301}-(2p_{211}+2q_{121})p_{201},\\
B_5 =& -2(p_{211}+q_{121})p_{301},\quad B_6 = \frac{1}{3}p_{122}+q_{032}-\frac{1}{3}(p_{111}+2q_{021})p_{021}, \\[1ex]
B_7 =&\,  (p_{111}+2q_{021})q_{031}-\frac{2}{3}(p_{211}+q_{121})p_{021},\quad B_8 =2(p_{211}+q_{121})q_{031}.
\end{split}
\end{equation}

Using $\rm(ii)$ and $\rm(iii)$ of Lemma \ref{lm1}, we get
\begin{equation}\label{b19}
\begin{split}
I^+_{5,1}=&\,hI^+_{1,1}+\frac{13}{12}I^+_{4,1}=\frac{13}{12}hI^+_{0,1}+hI^+_{1,1}+\frac{143}{126}I^+_{2,1},\\[1ex]
I^+_{1,3}=&\,\frac{3}{2}I^+_{5,1}-\frac{3}{2}I^+_{4,1}=\frac{1}{14}hI^+_{0,1}+\frac{3}{2}hI^+_{1,1}+\frac{11}{84}I^+_{2,1},\\[1ex]
I^+_{6,1}=&\,\frac{4}{3}hI^+_{2,1}+\frac{10}{9}I^+_{5,1}=\frac{130}{189}hI^+_{0,1}+\frac{10}{9}hI^+_{1,1}+(\frac{4}{3}h+\frac{715}{567})I^+_{2,1},\\[1ex]
I^+_{2,3}=&\,I^+_{6,1}-I^+_{5,1}=\frac{13}{189}hI^+_{0,1}+\frac{1}{9}hI^+_{1,1}+(\frac{4}{3}h+\frac{143}{1134})I^+_{2,1}.
\end{split}
\end{equation}
Then substituting \eqref{b14} and \eqref{b19} into \eqref{b18} yields
\begin{equation}\label{b19b}
M^+_2(h)=(B_1+C_3h)\,I^+_{0,1}+(B_2+C_4h)\,I^+_{1,1}+(C_5+C_6h)\,I^+_{2,1},
\end{equation}
where
\begin{equation}\label{b19a}
\begin{split}
C_3 =&\,\frac{4}{7}B_5+\frac{12}{7}B_6+\frac{1}{14}B_7+\frac{13}{189}B_8,\qquad C_4 =\frac{3}{2}B_7+\frac{1}{9}B_8,\\[0.5ex]
C_5 =&\, B_3+B_4+\frac{22}{21}B_5+\frac{1}{7}B_6+\frac{11}{84}B_7+\frac{143}{1134}B_8,\quad C_6 = \frac{4}{3}B_8.
\end{split}
\end{equation}
Then by \eqref{b13a} we can get
\begin{equation*}
M^+_2(h)=\bar c_0+b^+_{0}\bar c_1h^{\frac{5}{6}}+\bar c_2h
+b^+_{1}\bar c_3h^{\frac{7}{6}}+b^+_{0}\bar c_4h^{\frac{11}{6}}+\bar c_5h^2+\cdots,
\end{equation*}
where
\begin{equation}\label{b20}
\begin{split}
\bar c_0=&\,\frac{\sqrt{2}\pi}{19683}(2916B_1+2430B_2+2268B_3+2268B_4+2376B_5\\
&\qquad\quad  +324B_6+297B_7+286B_8),\\
\bar c_1=&-2B_1,\\
\bar c_2=&\frac{2\sqrt{2}\pi}{81}(81B_2+54B_3+54B_4+60B_5+18B_6+15B_7+14B_8),\\
\bar c_3=&\,B_1+2B_2,\\
\bar  c_4=&\,\frac{1}{44}(35B_1+42B_2+48B_3+48B_4-144B_6),\\
\bar c_5=&\,\sqrt{2}\pi(3B_7+2B_8).
\end{split}
\end{equation}
Note that $\bar c_j,\,j=0,1,\cdots,5,$ are linear in $p_{ij2}$ and $q_{ij2}$.
By \eqref{b18a} and \eqref{b20}, from $\bar c_0=\bar c_1=\bar c_2=\bar c_3=0$
we can get a unique solution
in $p_{122},\,q_{012},\,q_{112}$ and $q_{212}$, under which we have
\begin{equation}\label{b21}
\begin{split}
\bar c_4=&-\frac{5}{3}\,p_{021}(q_{121}+p_{211})\\[0.5ex]
&+\frac{1}{54}\,q_{031}(135p_{111}+244p_{211}+270q_{021}+244q_{121}),\\[0.5ex]
\bar c_5=
&\,\sqrt{2}\pi[-2\,p_{021}(q_{121}+p_{211})+\,q_{031}(3p_{111}+4p_{211}+6q_{021}+4q_{121})].
\end{split}
\end{equation}
Then if $q_{121}+p_{211}\neq0$, by \eqref{b21} solving $\bar c_4=0$ in $p_{021}$ yields
\begin{equation}\label{b22}
p_{021}=\frac{q_{031}(135p_{111}+244p_{211}+270q_{021}+244q_{121})}{90(q_{121}+p_{211})},
\end{equation}
and $\bar c_5=-\frac{64\sqrt{2}\pi}{45}q_{031}(q_{121}+p_{211})$.
Note that
\begin{equation*}
\mbox{det}\!\left(\frac{\partial(\bar c_0,\,\bar c_1,\,\bar c_2,\,\bar c_3,\,\bar c_4)}
{\partial (p_{122},\,q_{012},\,q_{112},\,q_{212},\,p_{021})}\right)\!\bigg{|}_{\eqref{b22}}
=\frac{2560\,\pi^2}{19683}(q_{121}+p_{211}).
\end{equation*}
Then by Lemma \ref{thm4} for $q_{031}(q_{121}+p_{211})\neq0$ system \eqref{sys} can have $8$ limit cycles
near $\Gamma_0$ with proper perturbations.

By \eqref{b19a} and \eqref{b20}, it is straightforward to get that
$\bar c_0=\bar c_1=\cdots=\bar c_5=0$ if and only if
\begin{equation}\label{b23}
B_1=B_2=B_7=B_8=0,\quad B_3=3B_6-B_4,\quad B_5=-3B_6,
\end{equation}
which is equivalent to  $B_1=B_2=0$, $C_3=C_4=C_5=C_6=0$.
Then by \eqref{b19b}, $M^+_2(h)\equiv0$ if and only if \eqref{b23} holds. By \eqref{b18a}, we can get all
the solutions \eqref{a4} for \eqref{b23}.
The proof is completed. \hfill $\square$

Next, we present the proof for Theorem \ref{thm3}.

%\begin{proof}
\noindent{\bf Proof of Theorem \ref{thm3}.}
To show the existence of ten limit cycles near $\Gamma_0$ in system \eqref{sys},
for the sake of simplicity
 we use the following perturbations
\begin{equation}\label{c1}
\begin{split}
P_1(x,y)=&\,y(p_{031}y^2+p_{211}x^2+p_{021}y+p_{111}x),\\[0.5ex]
P_2(x,y)=&\,xy(p_{021}(2q_{021}+p_{111})y+p_{212}x),\\[0.5ex]
 Q_1(x,y)=& -y^2(p_{211}x-q_{021}),\quad
 Q_2(x,y)= 0,%\\
\end{split}
\end{equation}
and $P_3(x,y)$ and $Q_3(x,y)$ are given in \eqref{b1}.
It is easy to get that $P_i$ and  $Q_i$, $i=1,2$ in \eqref{c1} satisfy \eqref{b17} and the condition $S_2$ of \eqref{a4}.
Then by Theorem \ref{thm1} we can get $M^\pm_k(h)\equiv0$, $k=1,2$.

It is straightforward to get
\begin{equation*}
\begin{split}
r_2=&\,\frac{2}{5}\,p_{031}(2q_{021}+p_{111})x^5-\frac{1}{2}p_{031}(2q_{021}+p_{111})x^4+\frac{1}{3}p_{211}(2q_{021}+p_{111})x^3\\
& +(p_{111}^2+3p_{111}q_{021}+2q_{021}^2-p_{212})x^2+p_{031}(2q_{021}+p_{111})xy^2,
\end{split}
\end{equation*}
which satisfies
\begin{equation*}
(Q_2+r_1Q_1)\mathrm{d}x-(P_2+r_1P_1)\mathrm{d}y=r_2\mathrm{d}H+\mathrm{d}{R_2},
\end{equation*}
where $r_1=-(p_{111}+2q_{021})x$ by \eqref{b17a}, and $R_2$ is a polynomial of degree $9$ in $(x,y)$.
Then using Fran\c coise's algorithm \cite{JP1996}, we obtain
\begin{equation}\label{c1a}
M^\pm_3(h)=\oint_{\Gamma_h^\pm}\widetilde{Q}_3\mathrm{d}x-\widetilde{P}_3\mathrm{d}y=\oint_{\Gamma_h^\pm}\left(\widetilde{Q}_3+\int \frac{\partial}{\partial x}\widetilde{P}_{3}\mathrm{d}y \right)\mathrm{d}x,
\end{equation}
where
\begin{equation*}
\begin{split}
\widetilde{Q}_3(x,y)=&\,Q_3(x,y)+r_1(x,y)Q_2(x,y)+r_2(x,y)Q_1(x,y),\\
\widetilde{P}_3(x,y)=&\,P_3(x,y)+r_1(x,y)P_2(x,y)+r_2(x,y)P_1(x,y).
\end{split}
\end{equation*}

Similarly as in the proof of Theorem \ref{thm1}, we shall use Lemmas \ref{lm1} and \ref{lm5} to study the asymptotic expansion of $M^+_3(h)$ in \eqref{c1a}.
By (ii) and (iii) of Lemma \ref{lm1}, we derive
\begin{equation}\label{c2}
\begin{split}
I^+_{3,3}=&\,\frac{12}{7}hI^+_{0,1}+\frac{1}{7}I^+_{2,1},\\[0.5ex]
I^+_{4,3}=&\,(\frac{442}{6237}h+\frac{48}{77}h^2)\,I^+_{0,1}+\frac{34}{297}hI^+_{1,1}
+(\frac{221}{1701}+\frac{988}{693}h)\,I^+_{2,1},\\[0.5ex]
I^+_{0,5}=&\, (-\frac{51250}{6237}h+\frac{240}{77}h^2)I^+_{0,1}
+\frac{170}{297}hI^+_{1,1}
-(\frac{110}{1701}-\frac{4940}{693}h)I^+_{2,1}.
\end{split}
\end{equation}
Then by \eqref{b14}, \eqref{b19} and \eqref{c2}, $M^+_3(h)$ can be simplified into the form
\begin{equation*}
M^+_3(h)=(B_1+B_2h+B_3h^2)I^+_{0,1}+(B_4+B_5h)I^+_{1,1}+(B_6+B_7h)I^+_{2,1},
\end{equation*}
where
\begin{equation}\label{c3}
\begin{split}
B_1 =& \,p_{103}+q_{013},\\
B_2 =& -\frac{1}{18711}(1782q_{021}^2+(102500p_{031}+891p_{111}-858p_{211})q_{021}\\
&+(51250p_{031}-429p_{211})p_{111}+891p_{212})p_{021}+\frac{4}{7}p_{123}+\frac{12}{7}q_{033},\\
B_3 =& \,\frac{80}{77}p_{031}(2q_{021}+p_{111})p_{021},\qquad B_4 = 2p_{203}+q_{113},\\
B_5 =& \,\frac{1}{891}(170p_{031}p_{111}+340p_{031}q_{021}+33p_{111}p_{211}-891p_{111}q_{021}\\
&+66p_{211}q_{021}-1782q_{021}^2-891p_{212})p_{021},\\
B_6 =& -\frac{1}{10206}(1782q_{021}^2+(440p_{031}+891p_{111}-858p_{211})q_{021}+(220p_{031}\\
&-429p_{211})p_{111}+891p_{212})p_{021}+q_{213}+\frac{1}{21}p_{123}+3p_{303}+\frac{1}{7}q_{033},\\
B_7 =& \,\frac{4}{2079}(2q_{021}+p_{111})(1235p_{031}+231p_{211})p_{021}.
\end{split}
\end{equation}
Then by \eqref{b13a} we get
\begin{equation*}
\begin{split}
M^+_3(h)=&\,c_0+b^+_{0}c_1h^{\frac{5}{6}}+c_2h+b^+_{1}c_3h^{\frac{7}{6}}+b^+_{0}c_4h^{\frac{11}{6}}
+c_5h^2+b^+_{1}c_6h^{\frac{13}{6}}+\cdots,
\end{split}
\end{equation*}
where
\begin{equation}\label{c4}
\begin{split}
c_0=&\,\frac{2\sqrt{2}\pi}{243}(18B_1 + 15B_4 + 14B_6),\qquad c_1=-2B_1,\\[0.5ex]
c_2=&\,\frac{2\sqrt{2}\pi}{243}(18B_2 + 243B_4 + 15B_5  + 162B_6 + 14B_7),\\[0.5ex]
c_3=&\,B_1+2B_4,\\[0.5ex]
c_4=&\,\frac{1}{44}(35B_1-88B_2 + 42B_4 + 48B_6),\\[0.5ex]
c_5=&\,\frac{2\sqrt{2}\pi}{3}(3B_5+2B_7),\\[0.5ex]
c_6=&-\frac{1}{208}(385B_1 - 208B_2 + 440B_4 - 416B_5 + 480B_6).
\end{split}
\end{equation}

By \eqref{c3} and \eqref{c4},  solving $c_0=c_1=\cdots=c_5=0$ in $p_{123}$, $q_{013}$, $q_{113}$, $q_{213}$,
$p_{212}$ and $p_{211}$ yields
\begin{equation}\label{c5}
\begin{split}
p_{123}=&\,\frac{4460}{891}p_{021}p_{031}p_{111}+\frac{8920}{891}p_{021}p_{031}q_{021}-3q_{033},\quad p_{211}=-\frac{155}{33}p_{031},\\ p_{212}=&\,\frac{85}{297}\,p_{031}p_{111}+\frac{170}{297}\,p_{031}q_{021}-p_{111}q_{021}-2q_{021}^2,\quad q_{013}=-p_{103},\\
q_{113}=&-2p_{203},\quad q_{213}=\frac{5}{891}p_{021}p_{031}p_{111}+\frac{10}{891}p_{021}p_{031}q_{021}-3p_{303},
\end{split}
\end{equation}
and then $c_6=-\frac{160}{297}\,(2q_{021}+p_{111})p_{031}p_{021}$,
\begin{equation*}
\mbox{det}\!\left(\frac{\partial(c_0,c_1,c_2,c_3,c_4,c_5)}
{\partial (p_{123},q_{013},q_{113},q_{213},p_{212},p_{211})}\right)\!\!\bigg{|}_{\eqref{c5}}\!\!\!\!
=-\frac{16384\sqrt{2}}{531441}\pi^3(2q_{021}+p_{111})p_{021}^2.
\end{equation*}
Then by Lemma \ref{thm4} system \eqref{sys} can have 10 limit cycles produced near $\Gamma_0$ with proper perturbations
for $|\varepsilon|$ sufficiently small and
parameters $p_{123}$, $q_{013}$, $q_{113}$, $q_{213}$, $p_{212}$, $p_{211}$ close enough to \eqref{c5} with $(2q_{021}+p_{111})p_{031}p_{021}\neq0$.
%\end{proof}
The proof is completed. \hfill $\square$

\section{Acknowledgments}
This work was supported by the National Natural Science Foundation of China (NSFC No. 11871042).

%\section*{References}

\end{document}